\documentclass[12pt]{amsart}
\usepackage{graphicx}

\newcommand{\Tet}{\mathrm{Tet}\begin{pmatrix}}

\newcommand{\cS}{\mathcal{S}}
\newtheorem{theorem}{Theorem}
\newtheorem{lemma}{Lemma}

\newtheorem{prop}{Proposition}
\newtheorem{definition}{Definition}

\newtheorem{remark}{Remark}  
\newtheorem{fact}{Fact}

\title{The Quantum Content of the Normal Surfaces in a Three-Manifold}

\author{Charles Frohman}

\address{Department of Mathematics, University of Iowa, Iowa City, IA
52242, USA}

\email{\tt frohman@math.uiowa.edu}

\author{Joanna Kania-Bartoszynska}

\address{Department of Mathematics, Boise State University, Boise, ID
83725, USA} 

\email{\tt kania@math.boisestate.edu}

\thanks{This material is based upon work supported by the National
Science Foundation under Grant No.0207030  and under Grant No.0204627}
\begin{document}
\begin{abstract}

The formula for the Turaev-Viro invariant of a $3$-manifold depends 
on a complex parameter $t$. When $t$ is not a root of unity, the formula 
becomes an infinite sum. This paper analyzes convergence of this sum  
when $t$ does not lie on the unit circle, in the presence of an efficient 
triangulation of the three-manifold. The terms of the sum can be indexed by
surfaces lying in the three-manifold. The contribution of a surface is 
largest when the surface is normal and when its genus is the lowest. 

\end{abstract}
\maketitle

\section{Introduction} This paper initiates the study of non-perturbative
quantum invariants of three-manifolds $M$  away from roots of unity. 
Turaev and Viro \cite{TV} defined  invariants of closed 3-manifolds
as  state sums depending on a complex parameter $t$. When $t$ is 
a root of unity this sum is finite.
At values of $t$ other than roots of unity the formula for
the Turaev-Viro invariant becomes an infinite sum. The partial sums
could oscillate wildly, so that even after
renormalizing the series does not converge. However, we are able to show 
that for a
special class of spines of some three-manifolds, the oscillation does not
occur and there is a limit. 

The key is
to see the invariant  as a sum over surfaces in the manifold.
An 
efficient ideal triangulation \cite{JR} of the manifold is one where
the only normal spheres and
tori are links of the boundary components of the manifold. Given an efficient
ideal triangulation we know what the invariant should be.
It is a sum over surfaces carried by a spine dual to the efficient
triangulation. 
The surfaces contribute to the sum in a way that fits the modern 
approach to
normal surface theory.
The study of normal surfaces \cite{JR} 
has been augmented by looking at surfaces that
aren't normal, and coming to an understanding of how a surface fails to be
normal. 
The farther a surface is from being normal, the less it contributes to
the sum for the invariant. The higher the genus of a surface the less it
contributes. In a very real sense, the Turaev-Viro invariant is a measure of
the normal surface theory of $M$.

In  section \ref{prel}  we broach preliminary concepts 
relating to special functions, and to spinal and normal surfaces in an
ideal triangulation of a $3$-manifold. This is followed by  section \ref{6jdet}
that studies  properties and limiting behavior 
of the  
$6j$-symbols. Section \ref{norandspi} is concerned with estimates  of
the contributions of spinal surfaces to the state sum.
The final section proves the result about the convergence of the
infinite state sum.

The authors would like to thank Ian Agol, William Jaco, Marc Lackenby, Sergei
Matveev, 
Dennis Roseman and Hyam Rubinstein for enlightening conversations about normal
surface theory.

\section{Preliminaries}\label{prel}
\subsection{Special Functions} The formulas in this section are taken
from \cite{KL}, however we use the variable $t$ instead of $A$.
Throughout this paper $t$ is a real number with $0<t<1$.
There are several functions of $t$  that we will work with. 
The first
is known as {\em quantized $n$},
\begin{equation} [n]=\frac{t^{2n}-t^{-2n}}{t^2-t^{-2}}.\end{equation}
The next is just a variation on the first,
\begin{equation}\Delta_n=(-1)^n\frac{t^{2n+2}-t^{-2n-2}}{t^2-t^{-2}}=(-1)^n[n+1].\end{equation}
There is quantized factorial, defined recursively by $[0]!=1$ and
\begin{equation} [n]!=[n][n-1]!.\end{equation}

\begin{figure}
\includegraphics{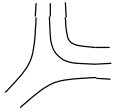}
\caption{Admissible triple $(2,3,3)$.}\label{adm}
\end{figure}

A triple of nonnegative integers $(a,b,c)$ is admissible if their sum
is even and they satisfy every possible triangle inequality.
Admissibility is the necessary and sufficient condition for the
existence of  a Kauffman triad on $a$, $b$ and $c$.
Suppose that $(a,b,c)$ is admissible.
Arrange $a$ points, $b$ points and
$c$ points on the sides of a triangle.
 There exists a system of disjoint
proper arcs  joining opposite sides of
the triangle having those points as their boundary.
Figure \ref{adm} shows the admissible triple
$(2,3,3)$. 
The number of strands running between the family of $a$ points, and the
family of $b$ points is $x_1=\frac{a+b-c}{2}$, between the $b$ points and the
$c$ points is $x_2=\frac{c+b-a}{2}$, and  between the $a$ points and the $c$
points is
$x_3=\frac{a+c-b}{2}$. Admissibility is equivalent to the statement that all
three functions are nonnegative and integral. We call $x_1$, $x_2$ and
$x_3$ the {\em strand numbers} of the triple $(a,b,c)$. 

You can always add two
admissible triples and the result will be admissible, but you cannot always
subtract them. However, 

\begin{prop}\label{subtract}
Suppose that $(a,b,c)$ and $(a',b',c')$ are admissible and
$x_i$ and $x'_i$ are the strand numbers corresponding to the two triples.
If $x_i\geq x'_i$ for all $i$, then  $(a-a',b-b',c-c')$ is
admissible. 
\end{prop}

\proof The strand numbers are linear functions of the triples, thus
the strand numbers for $(a-a',b-b',c-c')$ are integral and nonnegative.
\qed

If $(a,b,c)$ is an admissible triple, define, \cite{KL}
\begin{equation}
\theta(a,b,c)=(-1)^{\frac{a+b+c}{2}}\frac{[\frac{a+b+c}{2}+1]!
[\frac{a+b-c}{2}]![\frac{a+c-b}{2}]![\frac{b+c-a}{2}]!}{[a]![b]![c]!}.\end{equation}

\begin{figure}
\raisebox{-8pt}{\includegraphics{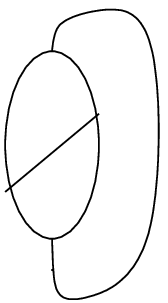}
\put(-50,16){$a$}
\put(-60,49){$b$}
\put(5,49){$e$}
\put(-23,70){$c$}
\put(-15,40){$d$}
\put(-40,49){$f$}}
\caption{Tetrahedral net}\label{tet}
\end{figure}

Suppose that a tetrahedral net has been labeled as in  Figure \ref{tet} 
where the letters are nonnegative integers and the triples appearing at each
vertex $v$ are admissible. 
The tetrahedral coefficient \cite{KL,MV} is the quantity
\begin{equation}\text{Tet}\begin{pmatrix} a & b & e \\ c & d &
f\end{pmatrix}=\end{equation}
\[\frac{\prod_v[x_{v,1}]![x_{v,2}]![x_{v,3}]!}{[a]![b]![c]![d]![e]![f]!}
\sum_{s=m}^M
\frac{(-1)^s[s+1]!}{[B_1-s]![B_2-s]![B_3-s]![s-A_1]![s-A_2]![s-A_3]![s-A_4]!},
\]
 where the $B_i$ are half the sums of the labels over the four cycles, the
 $A_i$ are half the sums of the labels at each vertex, $m$ is the maximum of
 the $A_i$ and $M$ is the minimum of the $B_i$.
The $x_{v,i}$ are the
 strand numbers of the admissible triple at the vertex $v$.

The {\em unitary $6j$ symbol} \cite{Ro} is the quantity:
\begin{equation} \left\{\begin{matrix} a & b & e \\ c & d &
f\end{matrix}\right\}_u=\frac{\text{Tet}\begin{pmatrix} a & b & e \\ c & d &
f\end{pmatrix}}
{\sqrt{\theta(a,d,e)\theta(b,c,e)\theta(a,b,f)\theta(c,d,f)}}.\end{equation}
A consequence of admissibility is that the denominator is a square root of a
positive number, so the formula is unambiguous.

Letting $q=t^4$, let $(x;q)_n=\prod_{i=1}^n(1-xq^{i-1})$.
We need the following fact \cite{GR}:
The function
\begin{equation}(x;q)_{\infty}=\prod_{i=1}^{\infty} (1-xq^{i-1})\end{equation} 
is well defined when $|q|<1$. In particular, $(q;q)_{\infty}$
is well defined.  
Notice that 
\begin{equation}\label{Delta} 
[n]=t^{-2n+2}\frac{1-q^{n-1}}{1-q},
\end{equation}
and
\begin{equation}
\Delta_n=(-1)^nt^{-2n}\frac{1-q^n}{1-q},  
\end{equation}
so that
\begin{equation}\label{fact}
[n]! = t^{-(n-1)n} \displaystyle\frac{(q;q)_n}{(1-q)^n},
\end{equation}
and
\begin{equation}\theta(a,b,c)= (-1)^{\frac{a+b+c}{2}}\frac{t^{-a-b-c}}{1-q}\  
\frac{(q;q)_{\frac{a+b+c}{2}+1}(q;q)_{\frac{a+b-c}{2}}
(q;q)_{\frac{b+c-a}{2}}(q;q)_{\frac{a+c-b}{2}}}
{(q;q)_a(q;q)_b(q;q)_c}.
\end{equation}

The quantities $\Delta_n$, $\theta(a,b,c)$ and 
$\text{Tet}\begin{pmatrix} a & b & e \\ c & d &
f\end{pmatrix}$ can be understood as the Kauffman brackets of colored graphs
\cite{KL}.

\subsection{Spines and Ideal Triangulations}
An {\em ideal triangulation}  \cite{BP} of the compact three-manifold $M$ is a 
union of  tetrahedra
joined along faces with their vertices removed so that the result is
homeomorphic to the interior of $M$.  

\begin{figure}
\includegraphics{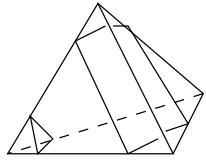}
\caption{Normal surface intersecting a tetrahedron}\label{norm1}
\end{figure}

\begin{figure}
\includegraphics{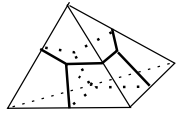}
\caption{Dual tetrahedral net}\label{net}
\end{figure}

A surface is {\em normal} with respect
to the triangulation \cite{JR} if it intersects each tetrahedron in 
triangles and quadrilaterals (quads) as in Figure \ref{norm1}.
Parameterize normal surfaces by their intersection
with the edges of the tetrahedra. Arrange these numbers in a $2\times
3$ array, so that  each column of the array is the number of points of 
intersection
of the normal surface with two opposite edges of the tetrahedron. More
specifically, there is
a tetrahedral net dual to the $1$-skeleton of the tetrahedron, lying on the
boundary of the tetrahedron, as pictured in Figure \ref{net}. 
The intersection of the normal surface with the boundary of the tetrahedron 
is a
family of circles carried by this net, the number of strands carried by an
edge of the net is the intersection number of the normal surface with the edge
of the tetrahedron transverse to the particular edge of the net. We form the
array of nonnegative integers just as if we were indexing a tetrahedral
coefficient, see Figure \ref{tet}, where the label on the edge is the number
of strands carried by that edge.  Let
$C_1,C_2,C_3$ be the sums of the columns of the array, named so
that
$C_1\geq C_2 \geq C_3$.

\begin{prop}\label{normalsurface}
 An array
of nonnegative integers corresponds to a normal surface if and only if
the integers assigned to the three edges around each face of the
tetrahedron form an admissible triple, and 
$C_1=C_2$. 
\end{prop}
\proof Think of a face of a
tetrahedron as a 
triangle, and the intersection of the normal surface with the triangle as a
system of arcs joining points on the three edges.
As we are joining the points on the three sides of the triangle by
nonintersecting arcs, the triple around each face must be admissible.
 The second condition follows from the fact that the intersection of
 a normal surface with any tetrahedron can
 only contain triangles and one type of quad. \qed

You cannot necessarily add normal surfaces, because if
two normal surfaces have different quads in the same tetrahedron, their double
curve sum may  no longer be normal. However, there is always a finite family of
normal surfaces so that every normal surface can be written as an integral sum
of those surfaces, with nonnegative coefficients.

Letting $\mathbb{N}$ denote the nonnegative integers, a {\em rational
cone} is the solution of
a family of linear homogeneous equations with integer coefficients
in $\mathbb{N}^k$ for some $k$.
An element of a rational cone is {\em irreducible} if it cannot be
written  as the sum of two elements of the
cone in a nontrivial way. It is a classical result that a rational
cone has only finitely many irreducible elements and they generate the cone additively. The set of
irreducibles is called a {\em Hilbert basis} for the cone.

The normal surfaces form a rational cone. The class of
normal surfaces is a subset of a more general class of surfaces, the
{\em spinal surfaces}.

Let \begin{equation}X=\{(x,y,z)\in \mathbb{R}^3| z=0 \ \text{or} \ (z\geq 0 \ \text{and} \
x=0) \ \text{or} \ (z\leq 0 \ \text{and} \ y=0) \} .\end{equation}
A subset $Y$ of the 3-manifold $M$ is {\em modeled} on $X$ if 
 for every 
$p \in Y$ there is an open neighborhood  $U$ of $p$, open set
$V\subset\mathbb{R}^3$   and
 a homeomorphism $\phi: U \rightarrow V$, with
$\phi^{-1}(X)=Y \cap U$. There is a decomposition  of $Y$ into 
vertices, (open) edges and (open) faces  coming from the natural 
decomposition of 
$X$ into vertices, edges and faces. 

We say $S \subset M$ is a {\em regular spine} \cite{Ma, Pi} if:
\begin{enumerate}
\item $S$ is  modeled on $X$.
\item $M -S $ is homeomorphic to $\partial M \times [0,1)$.
\item $S$ has at least one vertex.
\item Every edge of $S$ has a vertex in its closure.
\item Every face of $S$ is simply connected.
\end{enumerate}

\begin{figure}
\includegraphics{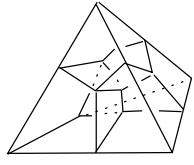}
\caption{Intersection of the spine with a tetrahedron}\label{spine}
\end{figure}

\begin{prop}
  Ideal triangulations and regular spines are in one to one correspondence up
  to isotopy via duality.
\end{prop}
\proof
For each ideal triangulation there is a regular spine. Put a vertex 
in the
center of each tetrahedron. Join the vertices in adjacent tetrahedra by edges,
and then form faces of the spine that intersect the edges of the triangulation
transversely and 
are bounded by the edges of the spine. 
The intersection of the spine with a tetrahedron is pictured in
Figure \ref{spine}.
Similarly for each regular spine there is an ideal triangulation, so that its
six edges intersect the six faces of the spine coming into the vertex 
transversely,
and each edge intersects exactly one face.\qed

Given a spine $S$ of $M$ and a  simple closed curve 
$\kappa \subset \partial M$ there is a possibly singular  annulus
$A_{\kappa} \subset M$ having $\kappa$ as one boundary component so that the
intersection of $A_{\kappa}$ 
with $S$ is the other boundary component of $A_{\kappa}$. 
The annulus is constructed by taking the closure of the points lying
over $\kappa$ in the 
product structure on $M-S$. The singularities of $A_{\kappa}$ come from the
fact that the map from $\partial M$ to $S$ given by following the lines of the
product structure is two to one along faces.
Since there is some ambiguity in the product structure we
can choose the annulus $A_{\kappa}$ so that it is in general position with
respect to the spine. This means its boundary  misses the vertices of the
spine, intersects the edges transversely and its only singular points are
transverse double points  occurring in the interior of the faces of the
spine.  If $\mathcal{C}$ is a system of disjoint  simple closed curves in 
$\partial M$
let $A_{\kappa}$, where $\kappa \in C$, be a system of disjoint annuli
corresponding to the curves in $C$ that is in general position with respect to
the spine $S$. 
The union $S(C)=S \cup \left( \cup_{\kappa \in C} A_{\kappa}\right)$ 
is called 
{\em the augmentation of the spine with respect to $C$}. Except for points on
$C$ the augmentation is still modeled on $X$, so it can be decomposed
into vertices, edges and faces just as a spine. If the spine is regular
then the faces of the augmentation are simply connected except for the
annular faces  with one boundary component a curve in $C$.

\begin{figure}
\includegraphics{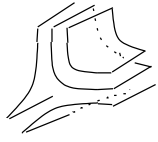}
\caption{Building a spinal surface}\label{glue}
\end{figure}

An {\em admissible coloring} of a spine is an assignment of a nonnegative 
integer to each face of the spine so that the integers assigned to the 
three faces meeting along each edge form an admissible triple. Given an
admissible coloring of the spine there is a {\em spinal surface} built as
follows. If the face $f$ carries the integer $u_f$ then take $u_f$ parallel
copies of $f$. Along the edges glue the faces together so that they look like
the Cartesian product of a triple of arcs at a vertex with an interval. The
triple $(2,3,3)$ occurring along an edge is shown in Figure \ref{glue}.
So far, the surface constructed intersects the boundary of a small ball at 
each vertex in a
collection of circles arranged along a tetrahedral net. To finish the
construction,  fill in the surface
inside each ball with a disk for each circle in the net.

Topologically, the spinal surfaces are those surfaces that intersect the
tetrahedra in disks, so that their intersection with any face of a tetrahedron
consists of arcs whose endpoints lie in distinct edges of the face.
The spinal surfaces form an additive cone, as the sum of two admissible
colorings is an admissible coloring. However, Euler characteristic is not
always additive under sum. Clearly, spinal surfaces are a larger class
than the normal surfaces associated to the dual triangulation.
We can identify the normal surfaces inside
the spinal surfaces by looking at the tetrahedral net at each vertex. In
specific at each vertex we can define the three column sums of the tetrahedral
net and order them so that $C_1 \geq C_2 \geq C_3$.
\begin{remark}
By  Proposition \ref{normalsurface} the surface is normal in the dual
ideal triangulation if and only if at each vertex
$C_1=C_2$. 
\end{remark}

The spinal surfaces form a rational cone. The proper domain
is the Cartesian product of copies of $\mathbb{N}$, one for each strand
number. The color on a face is the sum of the two adjacent strand
numbers along an edge of the face. The equations defining the cone
come from the 
requirement that the computed color of a face must be the same no
matter what edge of the face you compute it along. Thus we have the
following: 
\begin{fact}
There is a set of primitive spinal surfaces $\{F_i\}$ so that
every spinal surface can be written as a nonnegative  sum of the
$\{F_i\}$'s.
\end{fact}
Suppose now that $C$ is a system of disjoint simple closed curves in $\partial
M$ and $S(C)$ is an augmentation of the spine with respect to $C$. An
admissible coloring of $S(C)$ is defined the same way as an admissible coloring
of a spine except that the annular faces can only carry the color $1$. There
is once again a correspondence between admissible colorings and surfaces, but
now the surfaces have boundary equal to the union of the curves in $C$. 

In order to understand the Euler characteristic of a surface carried by a
spine or an augmented spine we need to understand how many circles there are
in a colored tetrahedral net.

\begin{prop}
  Suppose that a tetrahedral net has column sums $C_1 \geq C_2 \geq C_3$. 
The number of circles in the net is $\gcd{(C_1-C_2,C_1-C_3)}/2
+C_2+C_3-C_1 $. 
\end{prop}





\proof Unless
a tetrahedral net is of the form 
\begin{equation} 
\begin{pmatrix}
  a & b & a+b \\ a & b & a+b 
\end{pmatrix}\end{equation}
with $a$ and $b$ nonzero
then there is always a simple closed curve in the net that is the boundary of
one of the faces of the tetrahedron. 


Removing a  curve that bounds a face  does not change
$C_1-C_2$ or $C_1-C_3$, but $C_2+C_3-C_1$ is reduced by one. Remove
such curves until there are no more curves that bound faces. 
The remaining net will be of the form above.
If $\gcd{(C_1-C_2,C_1-C_3)}=2$ then the system consists of a single
curve. More generally, the number of components
is $\gcd{(a,b)}=\frac{\gcd{(C_1-C_2,C_1-C_3)}}{2}$. \qed

From this proof we see that the net is made up of circles that are
boundaries of faces along with multiple copies of a single type of circle that
appears in a tetrahedral net of type
\begin{equation} \begin{pmatrix}
  a & b & a+b \\ a & b & a+b 
\end{pmatrix}\end{equation}
where $a$ and $b$ are relatively prime. Alternatively, the boundary of a
simplex with its vertices removed is a four times punctured sphere. Any simple
closed curve is either boundary parallel or separates the surface into two
pairs of pants. The dearth of disjoint systems of simple closed curves
on a pair of pants causes all curves that are not triangles to be
parallel. The $a$ and $b$ can be understood in terms of geometric intersection
numbers with crosscuts.
Name such curves by the pair
$(a,b)$ where $a$ and $b$ are relatively prime and $a\leq b$.
 For each such
pair there are six or three different ways (depending on the
symmetries of the particular curve type) of labeling the tetrahedron
corresponding to the curve of type $(a,b)$.
We say that 
two $(a,b)$ curves are
{\em non-conflicting} if the curves are parallel in a regular neighborhood of
the $1$-skeleton of the tetrahedral net. 

\begin{prop} Euler 
characteristic of
spinal surfaces is additive when the two surfaces have the same $(a,b)$ types
at each vertex and those types are non-conflicting.
\end{prop}
\proof The surface that corresponds to the sum of the colorings is the
disjoint union of the surfaces corresponding to the two colorings. \qed

The type $(0,1)$ is a quad, the type
$(1,1)$ corresponds to an almost normal
surface \cite{JR}. Further types wind more and more around the 
tetrahedral net before closing up.




\section{$6j$-symbol Details}\label{6jdet}
\subsection{Bounding the $6j$-symbols.}

We begin with a universal bound on the size of the unitary
$6j$--symbols in terms of their entries.
\begin{prop}\label{6jest}
Let $C_1\geq C_2\geq C_3$ be the column sums of the unitary $6j$--symbol
$\left\{\begin{matrix} a & b & e \\ c & d &f\end{matrix}\right\}_u$ and
assume that $0<t<1$.  There exists a function $K(t)>0$ such that 
\begin{equation}\label{fundest}
\left|\left\{\begin{matrix} a&b&e \\ c&d&f\end{matrix}\right\}_u\right|
\leq K(t)t^{\frac{1}{2}(C_1-C_2)(C_1-C_3)+C_1}.
\end{equation}
\end{prop} 
\proof
After collecting and canceling terms,
$\left|\left\{\begin{matrix} a&b&e \\ c&d&f\end{matrix}\right\}_u\right|$
is equal to 
\begin{equation}
\frac{\sqrt{\prod_v\left|[x_{v,1}]![x_{v,2}]![x_{v,3}]!\right|}}
{\sqrt{[\frac{a+d+e}{2}+1]![\frac{b+c+e}{2}+1]![\frac{a+b+f}{2}+1]!
[\frac{c+d+f}{2}+1]!}} 
\left|\sum_{s= m}^{M}\frac{(-1)^s[s+1]!}
{\prod_{i=1}^3[B_i-s]!\prod_{j=1}^4[s-A_j]!}\right|.
\end{equation}
Using (\ref{fact}) this is further equal to
\begin{equation}\label{first}
t^p(1-q)^2
\sqrt{\frac{\prod_v (q;q)_{x_{v,1}} (q;q)_{x_{v,2}} (q;q)_{x_{v,3}}}
{(q;q)_{\frac{a+d+e}{2}+1} (q;q)_{\frac{b+c+e}{2}+1}
  (q;q)_{\frac{a+b+f}{2}+1} (q;q)_{\frac{c+d+f}{2}+1}}} \end{equation}
\[
\left|\sum_{s=m}^M t^{p_s} \frac{1}{1-q}
\frac{(-1)^s (q;q)_{s+1}}
{\prod_{i=1}^3(q;q)_{B_i-s}\prod_{j=1}^4(q;q)_{s-A_j} }  \right|.
\]
Here 
\begin{equation}\label{p}
 p=\frac{1}{2}(-a^2-b^2-c^2-d^2-e^2-f^2 +ae+ad+ab+af +be+bc+bf +ce+cd+cf + de+df)
\end{equation}
\[+a+b+c+d+e+f
\]
and 
\begin{equation}
p_s=6s^2+ a^2+b^2+c^2+d^2+e^2+f^2\end{equation} 
\[+ af+ac+ae+fc+fe+ce+be+bd+bf+ed+df+ab+ad+bc+cd\]
\[-2s(1+2a+2b+2c+2d+2e+2f).
\]
After completing the square, $p_s$ is
\begin{equation}\label{ps}
p_s=6\left(s-\frac{a+b+c+d+e+f+\frac{1}{2}}{3}\right)^2 
+\frac{1}{3}(a^2+b^2+c^2+d^2+e^2+f^2)\end{equation}
\[- \frac{1}{3}(af+ae+ac+ad+ab+fe+fc+ce+be+bd+bf+ed+df+bc+cd)
-\frac{2}{3}(a+b+c+d+e+f)-\frac{1}{6}.
\]
Combining all the factors of $(1-q)$ outside and the powers of $t$ inside
the sum,
formula (\ref{first}) can be simplified to
\begin{equation}\label{next}
(1-q)
\sqrt{\frac{\prod_v (q;q)_{x_{v,1}} (q;q)_{x_{v,2}} (q;q)_{x_{v,3}}}
{(q;q)_{\frac{a+d+e}{2}+1} (q;q)_{\frac{b+c+e}{2}+1}
  (q;q)_{\frac{a+b+f}{2}+1} (q;q)_{\frac{c+d+f}{2}+1}}} 
\end{equation}
\[
\left|\sum_{s=m}^M t^{p'_s}
\frac{(-1)^s (q;q)_{s+1}}
{\prod_{i=1}^3(q;q)_{B_i-s}\prod_{j=1}^4(q;q)_{s-A_j} }  \right|.
\]
where
\begin{equation}\label{power}
p'_s=6\left(s-\frac{a+b+c+d+e+f+\frac{1}{2}}{3}\right)^2 
-\frac{1}{6}(a^2+b^2+c^2+d^2+e^2+f^2)\end{equation}
\[
-\frac{1}{3}(ac+fe+bd)+\frac{1}{6}(ad+ae+de+bc+be+ce+ab+af+bf+cd+cf+df)
\]
\[+\frac{1}{3}(a+b+c+d+e+f)-\frac{1}{6}.
\]
The formula (\ref{next})  is the absolute value of an alternating
sum from $s=m$ to $s=M$. 
Take the quotient whose numerator is the summand
at $s+1$ and whose denominator is the summand at $s$, the result is,
\begin{equation} (-1)t^{12s-4(a+b+c+d+e+f)+4}
\frac{(1-q^{s+2})\prod_{i=1}^3(1-q^{B_i-s})}{\prod_{j=1}^4(1-q^{s+1-A_j}) }.\end{equation}
Take the absolute value, with
the effect of removing the (-1). In order to see that each one of these
quotients is smaller than the last, take the logarithm of the
result, giving:
\begin{equation}\label{loged}
 (12s-4(a+b+c+d+e+f)+4)\log{(t)}+\log(1-q^{s+2})
+\sum_{i=1}^3 \log{(1-q^{B_i-s})}
-\sum_{j=1}^4 \log{(1-q^{s+1-A_j})}
\end{equation}
Apply the Taylor series for $\log{(1-x)}$ to get
\begin{equation}\label{log}
 (12s-4(a+b+c+d+e+f)+4)\log{(t)}+\sum_{n=1}^{\infty}\frac{1}{n}\left(
-q^{n(s+2)}-\sum_{i=1}^3q^{n(B_i-s)}+\sum_{j=1}^4 q^{n(s+1-A_j)}
\right).
\end{equation}
As $t<1$ the first term gets smaller as $s$ increases. We analyze 
the sum over $n$ in (\ref{log}) term by term. As $n$ gets larger, $q^{n(s+2)}$ gets
smaller so $-q^{n(s+2)}$ gets larger. However, for each $i$ and $j$, 
$-q^{n(B_i-s)}$
and $+q^{n(s+1-A_j)}$ get smaller as $s$ increases. Furthermore the powers of
  $q$ appearing in any $q^{n(s+1-A_j)}$ are smaller than in $-q^{n(s+2)}$
    which means that the amount any one of them is decreasing is greater than
    the amount that $-q^{n(s+2)}$ is increasing, 
so each term is getting smaller. Therefore
    the sum over all $n$ is getting smaller and the quotients are decreasing.
Thus, the absolute value of the sum has a
unique maximum. 
Since the sum is alternating we  conclude that the absolute value of the
summand is less than or
equal to the largest term.
For any $n$, using $(q;q)_n\leq 1$ in the numerator and
$(q;q)_n\geq (q;q)_{\infty}$ in the denominator of (\ref{next})
together with 
the fact that the power of $t$ is the
largest when the exponent is the smallest, the  expression in equation
(\ref{next}) is smaller than 
\begin{equation}
(1-q)
\sqrt{\frac{1}{\left((q;q)_{\infty}\right)^4}}
t^{\mathrm{min}}
\frac{1}{\left((q;q)_{\infty}\right)^7},
\end{equation}
where ${\mathrm{min}}$ is the smallest value of $p'_s$.
Analyzing (\ref{power}) we can see that $p'_s$ is minimal when
 $s=M=(C_2+C_3)/2$.
Substituting we see that  
\begin{equation}
p'_s\geq \frac{(C_1-C_2)(C_1-C_3)}{2}+C_1.
\end{equation}
The final estimate is
\begin{equation} 
\left|\left\{\begin{matrix} a&b&e \\ c&d&f\end{matrix}\right\}_u\right|
\leq
\frac{1-q}{\left((q;q)_{\infty}\right)^9}\cdot
t^{\frac{(C_1-C_2)(C_1-C_3)}{2}+C_1}. 
\end{equation}
\qed

\subsection{ Some  important limits}

Let $\begin{pmatrix} a & b & e \\ c & d & f \end{pmatrix}$ be an admissible
labeling of the edges of a tetrahedron. For any nonnegative integer $k$, the
labelings 
$\begin{pmatrix} a+2k & b+2k & e+2k \\ c+2k & d+2k & f+2k
\end{pmatrix}$, 
  $\begin{pmatrix} a+2k & b+2k & e \\ c+2k & d+2k & f+2k
\end{pmatrix}$, and
$\begin{pmatrix} a+2k & b+2k & e \\ c+2k & d+2k & f
\end{pmatrix}$ are admissible.

\begin{prop}\label{lim}
  Given an admissible labeling
$\begin{pmatrix} a & b & e \\ c & d & f \end{pmatrix}$
 of a tetrahedral net,
the sequences
\begin{equation}\label{seq1}
 t^{-4k}\left\{ \begin{matrix} a+ 2k & b + 2k & e +2k \\ c + 2k & d +2k & f
      +2k \end{matrix} \right\}_u, 
\end{equation}
\begin{equation}\label{seq2} 
t^{-4k}\left\{ \begin{matrix} a+ 2k & b + 2k & e  \\ c + 2k & d +2k & f
      +2k \end{matrix} \right\}_u, 
\end{equation}
and
\begin{equation}\label{seq3}
 t^{-4k}\left\{ \begin{matrix} a+ 2k & b + 2k & e  \\ c + 2k & d +2k & f
       \end{matrix} \right\}_u 
\end{equation}
are  convergent. 
\end{prop}

We will only prove the first limit exists, the other two are similar.
The proof is based on the following elementary lemma.

\begin{lemma}\label{altsum}
  Suppose that $w(k)_n$ is a sequence of sequences so that for each fixed $k$,
  the sequence is alternating and converges to zero, and for fixed $n$ the
  sequence is convergent. Suppose further that there exists $N$ so
  that, independent of $k$, if $n\geq N$ then $|w(k)_n|\geq |w(k)_{n+1}|$.
 The sequence  \begin{equation}w(k)_{\infty}=\sum_{n} w(k)_n\end{equation} 
(depending on $k$) is convergent.
\end{lemma}

\proof This is an application of the proof of the alternating series test.
\qed

\proof (of  Proposition \ref{lim}) 
Recall Formula (\ref{next}). The strand numbers increase by $1$ each time $k$
increases by $1$ so the $(q;q)_{x_{v,i}}$ all converge to  
$(q;q)_{\infty}$ as $k$ goes to infinity. Similarly, the functions in the
denominator inside the radical all converge to $(q;q)_{\infty}$. Hence to 
prove the convergence we must only understand the quantities inside the sum.

Let $M(k)$, $m(k)$, $A_j(k)$, $B_i(k)$ and $p'(k)_s$ be the quantities
in (\ref{next})
associated to
\begin{equation}\left\{ \begin{matrix} a+ 2k & b + 2k & e +2k \\ 
c + 2k & d +2k & f +2k \end{matrix} \right\}_u,\end{equation}
as in the proof of Proposition \ref{6jest}. 
Let $n=M(k)-s$, and let
\begin{equation} w(k)_n=t^{p'(k)_s-4k}
\frac{(-1)^s (q;q)_{s+1}}
{\prod_{i=1}^3(q;q)_{B_i(k)-s}\prod_{j=1}^4(q;q)_{s-A_j(k)} } \end{equation}
for $n\leq M(k)-m(k)$, and $w(k)_n=0$ for $n>M(k)-m(k)$.

As $k$ increases by $1$, the $B_i(k)$ increase by $4$ and the $A_j(k)$ only
increase by $3$. So, $M(k)=M(0)+4k$, $m(k)=m(0)+3k$, $B_i(k)=B_i(0)+4k$, and 
$A_j(k)=A_j(0)+3k$. When $n=0$, $s=M(k)$ and 
\begin{equation}p'(k)_{M(k)}=\frac{(C(k)_1-C(k)_2)(C(k)_1-C(k)_3)}{2}
  +C(k)_1,\end{equation} 
which increases by $4$ when $k$ increases by $1$, so 
$t^{-4k+p'(k)_{M(k)}}$ is a constant.
We see that $w(k)_0$ is convergent. A
similar analysis
shows that for fixed $n$ the sequence $w(k)_n$ is convergent. The series is
clearly alternating.

We have already seen that for fixed $k$ the sequence $|w(k)_n|$ can have at 
most one
maximum, we just need to see that there is a bound on how big $n$ is at that
maximum  depending only on $t$, and 
$\begin{pmatrix} a & b & e \\ c & d & f \end{pmatrix}$. We do this by
looking at $\log{|w(k)_{n}/w(k)_{n+1}|}$ and seeing when it becomes
nonnegative. When $n>M(0)-m(0)+k$ then $w(k)_n=0$, so the maximum of
$|w(k)_n|$ has already occured. Hence we only need to understand the case when
the quotient $w(k)_{n}/w(k)_{n+1}$ is well defined.
To this end we substitute into Formula (\ref{loged}) to get,
\begin{equation}\label{logthing}
 \log{|w(k)_{n}/w(k)_{n+1}|}=
(12(M(0)-n)-4(a+b+c+d+e+f)+4)\log{(t)}\end{equation}
\[+\log(1-q^{M(0)+4k-n+2})
-\sum_{i=1}^4\log{(1-q^{M(0)+k-n+1-A_j(0)})}+\sum_{i=1}^3
\log{(1-q^{(B_i(0)-M(0)+n)})}.\]
Choose $N$ large enough so that, 
\begin{equation}\label{beegtrouble}
(12(M(0)-N)-4(a+b+c+d+e+f)+4)\log{(t)} >\end{equation}
\[-\log(1-q^{M(0)+4N-N+2})-\sum_{i=1}^3\log(1-q^{M(0)-N+2}).
\]
Notice that $-\sum_{i=1}^4\log{(1-q^{M(0)+k-n+1-A_j(0)})}>0$. 
Inequality (\ref{beegtrouble}) guarantees that the expression
(\ref{logthing}) is positive when 
$k=N$. 
Increasing k in $-\log(1-q^{M(0)+4k-N+2})$
 makes it smaller. Thus, by the argument from the proof of
 Proposition \ref{6jest}, the sequence $|w_k(n)|$ is monotone decreasing for
$n \geq N$. 

We have established the criterion for convergence from  Lemma \ref{altsum}.
\qed

Let the  limit of the  sequence (\ref{seq1}) from Proposition
\ref{lim} be denoted by: 
\begin{equation}\label{sixjinf}
\left\{\begin{matrix} a & b & e \\ c & d & f \end{matrix}\right\}_{\infty}
=\lim_{k\rightarrow \infty} t^{-4k}\left\{ \begin{matrix} a+ 2k & b +
    2k & e +2k \\ c + 2k & d +2k & f 
      +2k \end{matrix} \right\}_u.
\end{equation}
Similarly, denote the limits of the sequences (\ref{seq2}) and
(\ref{seq3}) by 
\begin{equation}\label{sixjinf1}
\left\{\begin{matrix} a & b & e_0 \\ c & d & f \end{matrix}\right\}_{\infty}
=\lim_{k\rightarrow \infty}
t^{-4k}\left\{ \begin{matrix} a+ 2k & b + 2k & e  \\ c + 2k & d +2k & f
      +2k \end{matrix} \right\}_u, 
\end{equation}
and
\begin{equation}\label{sixjinf2}
\left\{\begin{matrix} a & b & e_0 \\ c & d & f_0 \end{matrix}\right\}_{\infty}
=\lim_{k\rightarrow \infty}
 t^{-4k}\left\{ \begin{matrix} a+ 2k & b + 2k & e  \\ c + 2k & d +2k & f
       \end{matrix} \right\}_u. 
\end{equation}

\begin{remark}
\begin{equation}
 \left\{\begin{matrix} a & b & e \\ c & d & f \end{matrix}\right\}_{\infty}
= 
(1-q)(q;q)_{\infty} \sum_{u=0}^{\infty}(-1)^{\frac{C_2+C_3}{2}+u} 
\end{equation}
\[t^{6u^2+2(2C_1-C_2-C_3+1)u+\frac{(C_1-C_2)(C_1-C_3)}{2}+C_1}
\frac{1}{(q;q)_u(q;q)_{u+\frac{C_1-C_2}{2}}(q;q)_{u+\frac{C_1-C_3}{2}}}\]

The limits (\ref{sixjinf1}) and (\ref{sixjinf2}) are zero unless $a+c=b+d$.
\end{remark}

\section{Normal and Spinal Surfaces}\label{norandspi}
\subsection{ Analysis of the contribution of a surface}
For the remainder of this paper
$M$ will be a compact three-manifold with
non-empty connected  
boundary. Although the method works for a more general class of manifolds,
this assumption simplifies the arithmetic so that the ideas behind the
estimates are in the foreground.

\begin{definition} An ideal
triangulation $T$ whose only normal  spheres and  tori are the link of a
vertex is {\bf efficient}. An ideal triangulation
is {\bf $0$-efficient}
if and only if the only embedded, normal $2$-spheres are vertex
linking. 
\end{definition}
The $0$-efficient triangulations were studied in \cite{JR}. In particular,
it is shown there that any triangulation of a closed, orientable
irreducible $3$-manifold can be modified to a $0$-efficient
triangulation, or it can be shown that the $3$-manifold is one of
$S^3$, ${\mathbb R}P^3$ or $L(3,1)$. It is also shown that any
triangulation of a compact, orientable, irreducible and boundary
irreducible $3$-manifold with non-empty boundary can be modified to a
$0$-efficient triangulation. In the announced sequel to \cite{JR}
authors explore the concept of {\bf $1$-efficient manifolds}. They show that
the triangulations of irreducible, atoroidal, closed $3$-manifolds can
be obtained so that in addition to being $0$-efficient,  any
embedded normal torus is of a very special form or the $3$-manifold is
$S^3$, a lens space or a small Seifert fiber space.

Assume that $M$ has an efficient triangulation $T$.
Suppose that $S$ is the spine dual to $T$.  
Let $F$ be a surface carried by $S$. It is induced by an admissible
coloring of the 
spine. Let  $u_f$ denote the color assigned to the  
face $f$.
At each edge $e$ the three faces sharing that edge carry colors 
$a_e$, $b_e$ and $c_e$.
At each vertex there is a corresponding coloring of a tetrahedral net,
$ \left(\begin{matrix}a_v & b_v & e_v \\c_v & d_v & f_v
  \end{matrix}\right)$. 
Denote the column sums at vertex $v$ by $C_{1,v}\geq C_{2,v}\geq C_{3,v}$.

We can form the three strand numbers at each edge: 
$x_{e,1}=\frac{a_e+b_e-c_e}{2}$, $x_{e,2}=\frac{c_e+b_e-a_e}{2}$  and 
$x_{e,3}=\frac{a_e+c_e-b_e}{2}$. 
These are in fact linear functionals on the
space of spinal surfaces.  There is an arbitrariness to the choice of
which function is which, so fix this choice along each edge once and for
all. 
Similarly, at each vertex
we can form three linear functionals $S_{1,v}$, $S_{2,v}$ and $S_{3,v}$
corresponding to the column sums of the tetrahedral net at the
vertex. 
\begin{definition} If $C$ is a (possibly empty) set of simple closed
curves on the boundary of $M$, let $\mathcal{S}(C)$ denote the set
of spinal surfaces with respect to an augmentation of the spine
corresponding to $C$. 
For brevity, let $\mathcal{S}= \mathcal{S}(\emptyset)$.
A {\bf sector} $\mathcal{F}$ is determined  by fixing the
order of the values of the $S_{i,v}$ at each vertex (that is, deciding
which of the $S_{i,v}$'s is the largest column sum, $C_{1,v}$, etc.). 
\end{definition}  
Specifying these orderings  at all vertices breaks the space of 
spinal surfaces into $6^{\#v}$
sectors. 
Given any infinite sequence of spinal surfaces we can find a
subsequence that lives in one sector, because there are only finitely
many sectors.

\begin{prop}\label{subsurf}
  Suppose that the spinal surface $F$ lies in the sector $\mathcal{F}$.
Then every connected component of $F$ lies in the same sector.
\end{prop}

\proof Recall that the intersection of a spinal surface with a tetrahedron
consists of triangles along with one family of disks having a particular curve
type $(a,b)$. The triangles contribute the same to each column of the corresponding
symbol so any restriction on the sector comes from the curve type. Since all
components of $F$ are made up of a subset of the components of the
intersection of $F$ with each tetrahedron, they lie in any sector that
$F$ lies in (and maybe some other sectors too.) \qed 
\begin{definition}
A spinal surface $F$ is {\bf $k$-peelable} if $k$ is the maximum
non-negative integer such that $F$ can be written as
$F=F'+k\cdot\partial M$.
Use $\mathcal{S}_k(C)$
to denote the set of all surfaces in $\mathcal{S}(C)$ that are $k$-peelable.
Similarly, use $\mathcal{F}_k(C)$ to denote the $k$-peelable surfaces
in the sector $\mathcal{F}(C)$.
\end{definition}

There is a one-to-one correspondence between $\mathcal{S}_0$ and
$\mathcal{S}_{k}$ for any $k\geq 0$ given by $u_f\rightarrow u_f+2k$
for every $f$. Furthermore, this correspondence preserves sectors.

\begin{prop}
A spinal  surface is in $\cS_0$ 
if and only if some $x_{v,i}=0$. 
\end{prop}

\proof This follows from  Proposition \ref{subtract} (on being able
to subtract admissible triples). \qed 

Consequently, a spinal surface is $k$-peelable if and only if   
the minimum over all strand numbers $x_{v,i}$ is equal to $k$.

Let $Q:\mathcal{S} \rightarrow \mathbb{Z}$
be the function that assigns to each surface $F$,
\begin{equation}\label{Q} 
Q(F)= \sum_f -2u_f +
\sum_v \frac{1}{2}(C_{1,v}-C_{2,v})(C_{1,v}-C_{3,v})+C_{1,v}.
\end{equation}

\begin{prop} \label{cont}
\begin{itemize}
\item[{\bf(i)}] 
$\displaystyle{ {-2\chi(F)}\leq Q(F)}$
\item[{\bf(ii)}]  The function $Q(F)$ is super additive on any
sector. That is, for any surfaces $F$, $F'$ lying in the same sector, if
$C_{i,v}$ are the column sums corresponding to $F$ and $C'_{i,v}$ are the
column sums corresponding to $F'$ then
\begin{equation} Q(F+F')=Q(F)+Q(F')+
(C_{1,v}-C_{2,v})(C'_{1,v}-C'_{3,v})/2+(C'_{1,v}-C'_{2,v})
(C_{1,v}-C_{3,v})/2\end{equation}
\[\geq Q(F)+Q(F')\]
\item[{\bf(iii)}]  $Q(F)$ is bounded below on $\mathcal{S}_0$.
\item[{\bf(iv)}]  The level sets of $Q(F)$ on $\mathcal{S}_0$ are finite.
\item[{\bf(v)}]  The cardinality of the level sets of $Q$  on $\mathcal{S}_0$
  grows at most polynomially in the level.
\end{itemize}
\end{prop}

\proof
\begin{itemize}

\item[{\bf(i)}]  The Euler characteristic of the surface $F$ corresponding to the
coloring $u_f$ is
\begin{equation} \sum_f u_f -\sum_e \frac{a_e+b_e+c_e}{2}+\sum_v 
\gcd{(C_{1,v}-C_{2,v},C_{1,v}-C_{3,v})}/2
+C_{2,v}+C_{3,v}-C_{1,v} \end{equation}
Because each edge has exactly two ends we can redistribute the sum to
eliminate the sum over the edges. This yields,
\begin{equation}\label{eq1}  
\sum_f u_f+\sum_v 
\gcd{(C_{1,v}-C_{2,v},C_{1,v}-C_{3,v})}/2
+\frac{1}{2}C_{2,v}+\frac{1}{2}C_{3,v}-\frac{3}{2}C_{1,v}.
\end{equation}
Comparing (\ref{eq1}) to the right hand side of the inequality from
item {\bf(i)} 
we see that it is enough
to show that for each vertex $v$,
\begin{equation}\label{eq2}
-\gcd{(C_{1,v}-C_{2,v},C_{1,v}-C_{3,v})}-C_{2,v}-C_{3,v}+3C_{1,v} \leq 
\frac{1}{2}(C_{1,v}-C_{2,v})(C_{1,v}-C_{3,v})
+C_{1,v}.
\end{equation}
In the case that $C_{1,v}-C_{2,v}=0$ this reduces to 
$C_{1,v}\leq C_{1,v}$ thus the proposition is true. Assume that
$C_{1,v}-C_{2,v}>0$. The triples at each vertex are admissible so
$C_{1,v}-C_{2,v} \leq C_{1,v}-C_{3,v}$ are even and positive. Hence,
$\gcd{(C_{1,v}-C_{2,v},C_{1,v}-C_{3,v})}\geq 2$. Substituting this in
(\ref{eq2})
and putting everything on the right side, the inequality is equivalent to:
\begin{equation} \frac{1}{2}(C_{1,v}-C_{2,v}-2)(C_{1,v}-C_{3,v}-2)\geq 0.\end{equation}
Since we are assuming $C_{1,v}-C_{2,v}\geq 2$ and
$C_{1,v}-C_{3,v}\geq 2$ this is true. 

\item[{\bf(ii)}]  This is a direct computation from the formula. 

In what follows we would like to use this formula, to that end we write it
more compactly as follows. Letting
$F$ and $F'$ be surfaces in the same sector with
$\delta_v=C_{1,v}-C_{2,v}$, $\gamma_v=C_{1,v}-C_{3,v}$ being  associated with
$F$ and  $\delta_v'=C_{1,v}'-C_{2,v}'$, $\gamma_v'=C_{1,v}'-C_{3,v}'$ being associated
with $F'$, we have,
\begin{equation} Q(F+F')= Q(F)+Q(F')+\sum_v \frac{\delta_v\gamma_v'+\delta_v'\gamma_v}{2}.\end{equation}

\item[{\bf(iii)}] Since there are only
finitely many sectors, if  $Q$ is bounded below on each
sector, then it is bounded below on $\mathcal{S}_0$.
So assume we are working in a particular sector.
Suppose that $Q$ is not bounded below.
Starting with a surface with $Q<0$ we demonstrate the existence of
another surface of a particular
form with smaller $Q(F)$.
We then bound $Q$ below on surfaces of that form.

Suppose that $Q(F)<0$. Decompose $F$ as a union
$F_p$  of components with positive Euler characteristic 
and a union $F_n$ components with negative Euler
characteristic. 
Since $Q(F)<0$, the surface $F_p$ is nonempty.
By super-additivity we have that $Q(F)\geq
Q(F_p)+Q(F_n)$. Since $Q(F_n)\geq 0$ this implies that $Q(F_p)\leq
Q(F)$. Since $F_p$ is a subsurface of $F$, by Proposition
\ref{subsurf} it is in the same sector. 
So we can assume that we are working with a surface all of whose
components are spheres.

Next assume that $F$ has $\delta_v\geq 4$ for some $v$. Our estimate that
$Q(F)\geq -2\chi(F)$ tells us that if  $F$ has
a single component then $Q(F)\geq -4$. Using the fact that $\delta_v\geq 4$
for some vertex allows us to improve this to $Q(F)\geq -2$  Assume that
$F$ is not connected. We can 
then
write $F=F_1+F_2$ where the $F_i$ are from $\mathcal{S}_0$, and $F_2$ is
connected and has nonempty intersection with a small ball about $v$.  We use
$\delta_{v,1}$, $\delta_{v,2}$ to denote the differences between the largest 
column and second
largest column of these two surfaces at the vertex $v$, and $\gamma_{v,1}$ and
$\gamma_{v,2}$ to describe the difference between the largest column
and the smallest 
column. Note, 
$\delta_v=\delta_{v,1}+\delta_{v,2}$ and
$\gamma_v=\gamma_{v,1}+\gamma_{v,2}$. 
The super-additivity formula gives
\begin{equation}Q(F)=Q(F_1+F_2)=Q(F_1)+Q(F_2)+
\sum_v \frac{\delta_{v,1}\gamma_{v,2}+\delta_{v,2}\gamma_{v,1}}{2}.\end{equation}
Since $\delta_{v,1}+\delta_{v,2}\geq 4$ it follows that $Q(F_1)\leq Q(F)$ and 
it has
smaller $\delta_v$. Replace the surface $F$ with the surface $F_1$ and
continue until all $\delta_v\leq 2$.

Suppose  $F$ is a surface in $\mathcal{F}_0$  with all $\delta_v\leq 2$. Since
there are no normal spheres in $\mathcal{F}_0$ each sphere making up $F$ has
some $\delta_v=2$. Since $\delta_v$ is additive this means that there are no 
more
spheres in $F$ than there are vertices in the spine. Hence $Q$ is bounded
below by $-4(\text{\# vertices})$.

\item[{\bf(iv)}] It is enough to prove that the intersection of any level set
  with any sector is finite.
Suppose that $F_i$ is an infinite sequence of
spinal surfaces in a sector with $Q(F_i)=c$. If necessary we can pass to a
  subsequence so 
that  the strand numbers of the surfaces $F_i$ are
monotone increasing. There are two cases. 

{\bf Case 1} If the $C_{v,1}-C_{v,2}=\delta_v$ stay
bounded then we can further refine the sequence so that all these numbers are
constant. As the strand numbers are monotone increasing we can subtract the
first term of the sequence from every subsequent term to get a new
sequence of  spinal 
surfaces which are normal. The
values of $Q(F_i)$ are bounded  
below (by item {\bf(iii)}), hence there is an infinite sequence of surfaces
with 
the same Euler characteristic. Since these surfaces all  have some
strand number 
$0$, and the triangulation is efficient  they can be written as a sum of a 
finite list of normal surfaces so
that none of the surfaces has positive or zero Euler characteristic.
This is a contradiction, as
their Euler characteristic is increasing.

{\bf Case 2} If some $C_{v,1}-C_{v,2}=\delta_v$ is unbounded we refine
the sequence
so that the $\delta_v$ are monotone increasing and the strand numbers are 
monotone
increasing. Let $v$ be a vertex where  the $\delta_v$ are unbounded, and assume
that the first surface in the sequence has $\delta_v>0$. If not, just start
later. Subtracting the first surface from every surface in the
sequence
the super-additivity formula informs us that this is
a sequence of surfaces in $\mathcal{S}_0$ such that $Q$ is not bounded
below. This contradicts item $({\bf iii)}$.

\item[{\bf(v)}]  
If $V$ is a finite dimensional free $\mathbb{Z}$-module and $v_i$ is
a basis, we can define $N:V \rightarrow \mathbb{Z}$ by
\begin{equation} N(\sum_i c_i v_i)=\sum_i |c_i|.\end{equation}
The cardinality of the set of elements in $V$ with $N(v)\leq n$ is
less than or equal to a 
polynomial in $n$. Fixing a sector $\mathcal{F}$ there is a finite family of
surfaces $F_i$ that generate the surfaces in $\mathcal{F}$ as an integer cone.
As there are only finitely many surfaces $F$ with $Q(F)\leq0$, there is an
integer $K$ so that for any $\sum_i c_i F_i$, if some $c_i\geq K$ then
$Q(\sum_i c_i F_i)>0$. Let $S_j$ be the set of all surfaces $\sum_i c_i F_i$,
so that some $c_i$ is between  $K$ and $2K$ and the 
other $c_i$ are between $0$ and $K-1$. It is clear that all but finitely many
surfaces in
$\mathcal{S}$ can be written as a positive sum of these surfaces. Form
a free $\mathbb{Z}$-module with basis $v_j$
corresponding to the $S_j$ and define
a map from the nonnegative integer sums of the $v_j$ to $\mathcal{S}$ by
sending the $v_j$ to the $S_j$. This map is onto all but a finite subset
of $\mathcal{F}$. Also,
\begin{equation} N(\sum_i c_j v_j)\leq Q(\sum_j c_jS_j),\end{equation}
so the level set $Q(S)=n$ is the image of a subset of $V$ contained inside the
set $N(v)\leq n$. Therefore the level sets of $Q$ grow at most polynomially in $n$.

\end{itemize} \qed

Now suppose that $C$ is a system of simple closed curves on $\partial M$.
We consider colorings of an augmentation of the spine corresponding to
$C$. Let $\chi (f)$ denote the Euler characteristic of the face $f$.
Note that $\chi (f)=1$ if $f$ is an open disk, and  $\chi(f)=0$ if $f$ is
an annulus.

The space of surfaces corresponding to admissible colorings of the augmented
spine is much like the space of spinal surfaces, except you can't add two
augmented colorings. However, you can add an augmented coloring and any
coloring of the original spine. We can divide the space of colorings of the
augmented spine into sectors just like we did for spinal surfaces, and we can
define $k$-peelable. Let $\mathcal{F}(C)$ be the surfaces in a
sector coming from 
colorings of an augmentation of the spine,
and denote by $\mathcal{F}$ the corresponding sector in space of
spinal surfaces associated to the original spine. Use  
$\mathcal{F}(C)_k$ to denote the  $k$-peelable surfaces in that sector.
Define $Q$ from the space of surfaces
corresponding to admissible colorings of the augmented spine to the counting
numbers by,
\begin{equation}\label{newQ} 
Q(F)= \sum_f -2 \chi (f) u_f +
\sum_v \frac{1}{2}(C_{1,v}-C_{2,v})(C_{1,v}-C_{3,v})+C_{1,v}.
\end{equation}

\begin{prop}\label{aug}
\begin{itemize}
\item[{\bf(i)}]  $ {-2x_F}\leq Q(F)$
\item[{\bf(ii)}]  The function $Q(F)$ is super additive on sectors.
 If
$C_{i,v}$ are the column sums corresponding to $F \in \mathcal{F}(C)$
 and $C'_{i,v}$ are the 
column sums corresponding to $F'\in \mathcal{F}$ then
\begin{equation} Q(F+F')=Q(F)+Q(F')+
(C_{1,v}-C_{2,v})(C'_{1,v}-C'_{3,v})/2+(C'_{1,v}-C'_{2,v})(C_{1,v}-C_{3,v})/2\end{equation}  
\[\geq Q(F)+Q(F').\]
\item[{\bf(iii)}]  $Q(F)$ is bounded below on $\mathcal{S}(C)_k$.
\item[{\bf(iv)}]  The level sets of $Q(F)$ on $\mathcal{S}(C)_k$ are finite.
\item[{\bf(v)}]  The cardinality of the level sets of $Q$  on $\mathcal{S}(C)_k$ grows
at most polynomially in the level.
\end{itemize}
\end{prop}

\proof The proof is an extension of the proof of  Proposition \ref{cont}. The
first two parts follow directly from the formula for $Q$. 

The third part we
argue  as follows. First get the estimate on $\mathcal{S}(C)_0$
by working in sectors. 
Given a surface $F \in \mathcal{F}(C)_0$ it can be written as a sum
of a surface $F_1$ such that each of its components has nonempty boundary and
a surface $F_2$ each
component of which is closed. From Proposition \ref{cont} we have a lower
bound for $Q(F_2)$, from inequality {\bf(i)}  we
can bound $Q(F_1)$ below 
by $-2$ times the number of components in $C$. By super-additivity we have
bounded $Q$ from below on $\mathcal{F}(C)_0$. 

To bound $Q$ below on $\mathcal{S}(C)_k$ use the one-to-one correspondence
between surfaces in $\mathcal{S}(C)_k$ and surfaces in $\mathcal{S}(C)_0$
obtained by adding $k$ copies of $\partial M$. Once again we bound the value
of $Q$ on $k$ parallel copies of the boundary using the inequality from
item {\bf(i)} 
and then use the bound on $\mathcal{S}(C)_0$ and super-additivity
on sectors.

The proofs of items {\bf(iv)}  and {\bf(v)}  are completely analogous
to the proofs in  Proposition \ref{cont}.
\qed

\subsection{Summing Over  $k$-peelable Surfaces}
 Let $C$ be a system of simple closed curves in $\partial M$,  let
$S$ be a spine that is dual to an efficient triangulation of $M$
and let $S(C)$ be an augmentation of $S$ with respect to $C$. 
Given a coloring  $F$  of the augmented spine $S(C)$  let the $u_f$,
$a_e$, $b_e$, $c_e$, $a_v$, $b_v$, $c_v$, $d_v$, $e_v$, $f_v$ and $\chi(f)$ be
defined as before. 
Also, let $\chi(e)=1$ if the edge $e$ has some vertex in its closure and 
let $\chi(e)=0$ otherwise ($\chi(e)$ is the Euler characteristic of the
edge $e$).
\begin{definition}
The contribution of $F$ is defined to be
\begin{equation} E(F)=\frac{\prod_f \Delta_{u_f}^{\chi(f)} \prod_v \Tet a_v & b_v & e_v \\
c_v & d_v & f_v \end{pmatrix}}{\prod_e \theta(a_e,b_e,c_e)^{\chi(e)}}.\end{equation}
\end{definition}
Notice that faces 
and edges of the spine contribute to $E(F)$ unless they are annular
or belong to the simple closed curves on the boundary respectively. 
Each vertex is an endpoint of four edges and each edge
that counts in the contribution of a surface has two ends. We can thus
collect the tetrahedral coefficient at each vertex with the thetas to
reparse this product as
\begin{equation}\label{reparse} 
E(F)= \prod_f \Delta_{u_f}^{\chi(f)}\prod_v \left\{ 
    \begin{matrix}
     a_v & b_v & e_v \\
c_v & d_v & f_v 
    \end{matrix}\right\}_u.
\end{equation} 

There is a map $S(C)_0 \rightarrow S(C)_k$ that adds $k$ copies of the
boundary of $M$ (as a union of triangles near the vertex). This map is one to one
and onto. If the largest color corresponding to $F$ is $N$ then the largest
color corresponding to $F+k\partial M$ is $N+2k$.
We define 
\begin{equation}\label{Ek}
E_k(F)=E(F+k\partial M).
\end{equation}
Since
$\chi(M)=\#f-\#v$,
\begin{equation}\label{kpeel}Q(F+k\partial
M)=Q(F)-4k\chi(M).\end{equation}
Using results of Proposition \ref{lim} about limits of $6j$ symbols we have,
\begin{prop}\label{Es}
For every surface $F \in \mathcal{S}(C)_0$, the limit
\begin{equation}\label{Einf}
\lim_{k\rightarrow\infty}t^{4k\chi(M)}E_k(F)
\end{equation}
 exists. When $C=\emptyset$, it is equal to
\begin{equation}\label{limeinf}
E_{\infty}(F)= 
\prod_f (-1)^{u_f}
\frac{t^{-2u_f}}{1-q}
\prod_{v}\left\{  \begin{matrix}
     a_v & b_v & e_v \\
c_v & d_v & f_v 
    \end{matrix}\right\}_{\infty}. 
\end{equation}
\end{prop}

\proof

Assume first that $C=\emptyset$, thus $\chi(f)=1$ for all $f$.
Given a surface $F \in \mathcal{S}(C)_0$  and  $k>0$, 
use (\ref{reparse}) together with (\ref{Delta})
to express  
\begin{equation}\label{firste}
E_k(F)=\prod_{f} (-1)^{-u_f-2k}t^{-2u_f-4k} \frac{1-q^{u_f+2k}}{1-q}  
\prod_v t^{4k}t^{-4k}
\left\{ 
    \begin{matrix}
     a_v+2k & b_v+2k & e_v+2k \\
c_v+2k & d_v+2k & f_v+2k 
    \end{matrix}\right\}_u.
\end{equation}
Since $\chi(M)=\#f-\#v$, equation (\ref{firste}) can be rewritten as
\begin{equation}\label{seconde}
t^{-4k\chi(M)}\prod_f (-1)^{u_f}t^{-2u_f}\frac{1-q^{u_f+2k}}{1-q} 
\prod_v t^{-4k} 
\left\{  \begin{matrix}
     a_v+2k & b_v+2k & e_v+2k \\
c_v+2k & d_v+2k & f_v+2k 
    \end{matrix}\right\}_u.
\end{equation}
By Proposition \ref{lim}, along with the fact that
$\lim_{k\rightarrow\infty}\frac{1-q^{u_f+2k}}{1-q}=\frac{1}{1-q}$,  
limit  (\ref{Einf}) exists and is
given by the formula (\ref{limeinf}).

In the case when $C\neq\emptyset$ the argument is similar. The product
in (\ref{firste}) must be taken over all faces $f$ with $\chi(f)\neq
0$ and for some of the vertices $v$ we need to consider the limit of 
sequences (\ref{seq2}) or (\ref{seq3}) instead of the sequence
(\ref{seq1}) as in equation (\ref{seconde}).

\qed

Let $\mathcal{S}(C)_k^N$ be the subset of $\mathcal{S}(C)_k$
 where the largest color $u_f$ is less than or equal to $N$ and let
$\mathcal{S}(C)_k^{T(N)}$ be the subset of $\mathcal{S}(C)_k$
where the largest $u_f$ is greater than $N$, the {\em tail} of the set.
Clearly, 
\begin{equation}\mathcal{S}(C)_k=\mathcal{S}(C)_k^N\cup \mathcal{S}(C)_k^{T(N)}.\end{equation}

\begin{lemma}\label{tail}
  For every $\epsilon >0$ there is $N$ so that for all $k$,
\begin{equation} 
\sum_{F \in \mathcal{S}(C)_k^{T(N+2k)}} |E(F)|< t^{-4k \chi(M)} \epsilon,
\end{equation}
where $\chi(M)$ is the Euler characteristic of the manifold $M$.

Moreover, 
for every $i\geq 0$, the limit 
\begin{equation}\label{Fki}
\lim_{k\rightarrow\infty}t^{4k\chi(M)}\sum_{F \in
  \mathcal{S}(C)_k^{k+i}} |E(F)|
\end{equation}  
exists.

\end{lemma}

\proof 
Using (\ref{reparse})
 along with the estimate from Proposition
\ref{6jest}, we see that, 
\begin{equation} |E(F)| \leq D(t,M,C) t^{Q(F)},\end{equation}
where $D(t,M,C)$ is a number that only depends on $t$, the manifold $M$ and
the augmentation of the spine corresponding to $C$. From Proposition \ref{aug}
the function $Q(F)$ is bounded below by some $Q_0 \in \mathbb{Z}$,
and has finite level sets, so that the
level set where $Q$ takes on the value $n$ has its cardinality bounded above
by a polynomial $p(n)$. Comparing with 
\begin{equation} \sum_{n\geq Q_0} p(n)t^n,\end{equation}
the series 
\begin{equation}\sum_{F \in \mathcal{S}(C)_0} |E(F)|\end{equation}
is absolutely summable. This means that for each $\epsilon>0$ there is
$N$ so that 
\begin{equation} \sum_{F \in \mathcal{S}(C)_0^{T(N)}} D(t,M,C) t^{Q(F)}<  \epsilon.\end{equation}
Using equation (\ref{kpeel})
we have
\begin{equation} \sum_{F \in \mathcal{S}(C)_k^{T(N+2k)}} D(t,M,C) t^{Q(F)}<  
t^{-4k\chi(M)}\epsilon.\end{equation} 

Combining the above argument with Proposition \ref{Es} yields the
existence of the limit (\ref{Fki}).
\qed

\begin{remark}\label{restatedlemma}
  The first part of this lemma can be restated as follows: for every 
$\epsilon>0$ there exists $N$ so that independent of $k$,
\begin{equation} \sum_{F \in \mathcal{S}(C)_0^{T(N)}} 
t^{4k\chi(\partial M)} |E_k(F)|<\epsilon.\end{equation}

\end{remark}

\begin{prop}\label{THEZ}
Let 
\begin{equation}
Z_k(M)=\sum_{F \in \mathcal{S}(C)_k} E(F)=
\sum_{F \in \mathcal{S}(C)_0} E_k(F),
\end{equation}
and
\begin{equation}
|Z|_k(M)=\sum_{F \in \mathcal{S}(C)_k} |E(F)|=
\sum_{F \in \mathcal{S}(C)_0} |E_k(F)|=.
\end{equation}
For each $k$, $Z_k(M)$ and $|Z|_k(M)$ are well defined.
Moreover, the limits 
$Z_{\infty}(M)=\lim_{k\rightarrow \infty} t^{4k\chi(M)}Z_k(M)$
and $|Z|_{\infty}(M)=\lim_{k\rightarrow \infty} t^{4k\chi(M)}|Z|_k(M)$
exist.
\end{prop}

\proof
The well defined part of the proposition follows directly from Lemma
\ref{tail}.

In order to prove convergence of $Z_k(M)$,
choose $\epsilon >0$. There exists $N$ so that for all $k$
\begin{equation} \sum_{F \in \mathcal{S}(C)_k^{T(N+2k)}} D(t,M,C) t^{Q(F)}<  
t^{-4k\chi(M)}\epsilon/4.\end{equation}
By Proposition \ref{Es} there is a $K$ so that if $k_1,k_2>K$ then
\begin{equation} | t^{4k_1\chi(M)}\sum_{F \in\mathcal{S}(C)_{k_1}}^NE(F)- 
t^{4k_2\chi(M)}\sum_{F \in \mathcal{S}(C)_{k_2}}^NE(F)|<\epsilon/2.\end{equation}
This means that
\begin{equation} |t^{4k_1\chi(M)}Z_{k_1}(M)-t^{4k_2\chi(M)}Z_{k_2}(M)|\leq \end{equation}
\[ | t^{4k_1\chi(M)}\sum_{F \in \mathcal{S}(C)_{k_1}}^NE(F)- 
t^{4k_2\chi(M)}\sum_{F \in\mathcal{S}(C)_{k_2}}^NE(F)|+\]
\[|t^{4k_1\chi(M)}\sum_{F \in \mathcal{S}(C)_{k_1}^{T(N+2k_1)}}E(F)|+
|t^{4k_2\chi(M)}\sum_{F \in \mathcal{S}(C)_{k_2}^{T(N+2k_2)}}E(F)|\leq\]
\[ \epsilon/2 +\epsilon/4+\epsilon/4.\] As the sequence is Cauchy it
converges. The same proof works for $|Z|_{\infty}$. \qed

\section{The Invariant Sums}
In the section we analyze the sum of contributions of all spinal
surfaces in the three-manifold $M$ with an efficient triangulation.

Given any integer $r\geq 3$,  all the
special functions, $\Delta_n$, $\theta(a,b,c)$,
$\text{Tet}\begin{pmatrix} a & b & e \\ c & d &
f\end{pmatrix}$, are well defined for $t=e^{\frac{\pi i}{2r}}$
whenever $a,b,c,d,e,f\leq r-1$ and the 
condition $a+b+c\leq 2r-4$ is added to the definition of
admissibility. Given a system $C$ of disjoint simple closed curves in
$\partial M$ and  an augmentation $S(C)$ of the spine dual to the
triangulation of $M$,
the (finite) sum over all $r$-admissible colorings of the faces
of $S(C)$,
\begin{equation}\label{turvir}
\sum_{\text{$r$-admissible colorings of $S(C)$}}
  \frac{\prod_ f\Delta_{u_f}^{\chi(f)} \prod_v \Tet a_v & b_v & e_v \\
    c_v & d_v & f_v 
\end{pmatrix} } {\prod_e \theta(a_e,b_e,c_e)^{\chi(e)}},
\end{equation}
is a coefficient (corresponding to $C$) of a vector-valued invariant
associated to $M$ by the topological quantum field theory underlying
the Turaev-Viro invariant of $M$ at level $r$. Our idea is the extend
the invariant away from the roots of unity. The first major step is to
analyze the convergence of the infinite sums like (\ref{turvir}),
where $t=e^{\frac{\pi i}{2r}}$ is replaced by any $0<t<1$ and the
colorings are admissible.

\begin{theorem}
Let $S(C)^N$ denote the set of admissible colorings of $S(C)$
with all $u_f\leq N$.
\begin{itemize}
\item[{\bf(i)}]   If the Euler characteristic of $M$ is negative then
\begin{equation}\label{bigen}
 \sum_{\text{admissible colorings $u_f$ of $S(C)$}}
  \frac{\prod_ f\Delta_{u_f}^{\chi(f)} \prod_v \Tet a_v & b_v & e_v \\
    c_v & d_v & f_v 
\end{pmatrix} } {\prod_e \theta(a_e,b_e,c_e)^{\chi(e)}}
\end{equation}
converges absolutely.
\item[{\bf(ii)}] If $\chi(M)=0$ then
\begin{equation}\label{knotcomp} 
\lim_{N\rightarrow \infty}\frac{1}{N}
\sum_{S(C)^N}
  \frac{\prod_{f}\Delta_{u_f}^{\chi(f)} \prod_v \Tet a_v & b_v & e_v
    \\ c_v & d_v & f_v 
\end{pmatrix}}{\prod_e \theta(a_e,b_e,c_e)^{\chi(e)}}
\end{equation}
exists and is equal to $Z_{\infty}(M)=\sum_{F\in S_0(C)}E_{\infty}(F)$
which converges  absolutely. 
\item[{\bf(iii)}] 
If $\chi(M)=1$ then
\begin{equation}\label{sphere} 
\lim_{N\rightarrow \infty} t^{8N}
\sum_{S(C)^{2N}}
  \frac{\prod_{f}\Delta_{u_f}^{\chi(f)} \prod_v \Tet a_v & b_v & e_v
    \\ c_v & d_v & f_v 
\end{pmatrix}}{\prod_e \theta(a_e,b_e,c_e)^{\chi(e)}}
\end{equation}
exists. 
Given a spinal surface $F$, let $m(F)$ denote the least even number
greater than or equal to the  maximal
color corresponding to $F$.
The limit (\ref{sphere}) is equal to the sum of the absolutely
convergent series: 
\begin{equation}\label{spherelim}
\frac{1}{1-q} 
\sum_{F\in S_0(C)}
t^{4m(F)}
E_{\infty}(F).
\end{equation} 
\end{itemize}
\end{theorem}

\proof 
\begin{itemize}
\item[{\bf(i)}]
We need to show that the sequence of partial sums of the absolute
values of the  series
(\ref{bigen}) converges, that is,  
\begin{equation} \lim_{N\rightarrow \infty}
\sum_{S(C)^N}
 \left| \frac{\prod_{f}\Delta_{u_f}^{\chi(f)} \prod_v \Tet a_v & b_v & e_v
    \\ c_v & d_v & f_v 
\end{pmatrix}}{\prod_e \theta(a_e,b_e,c_e)^{\chi(e)}}\right|\end{equation}
exists.
Notice that 
\begin{equation}\sum_{S(C)^N}
\left|  \frac{\prod_{f}\Delta_{u_f}^{\chi(f)} \prod_v \Tet a_v & b_v & e_v
    \\ c_v & d_v & f_v 
\end{pmatrix}}{\prod_e \theta(a_e,b_e,c_e)^{\chi(e)}}\right|=\end{equation}
\[\sum_k 
\sum_{F \in \mathcal{S}(C)_k^N} |E(F)|<
\sum_k |Z|_k(M).\]
Proposition \ref{THEZ} implies that the series $\sum_k|Z|_k(M)$
converges by comparison with the series $\sum_kt^{-4k\chi(M)}$.

\item[{\bf(ii)}]
First, regroup the finite sum in (\ref{knotcomp}) according to
$k$-peelable surfaces. That 
is, use the fact that $S(C)^N$ is a disjoint union of subsets
$S(C)_k^N$ with $k=0,\dots ,N$ (since $\mathcal{S}(C)_k^N$ is empty for
$k>N$). Thus,
\begin{equation}\label{regroup}
\frac{1}{N}
\sum_{S(C)^N}
  \frac{\prod_{f}\Delta_{u_f}^{\chi(f)} \prod_v \Tet a_v & b_v & e_v
    \\ c_v & d_v & f_v 
\end{pmatrix}}{\prod_e \theta(a_e,b_e,c_e)^{\chi(e)}}=
\frac{1}{N}
\sum_{k=0}^{N}\sum_{F \in \mathcal{S}(C)_k^N} E(F).
\end{equation}
By Proposition \ref{THEZ} we can find $K$ so that for all $k>K$ we
have 
\begin{equation}\label{zbound}
|Z_k(M)-Z_{\infty}(M)|<\frac{\epsilon}{4}.\end{equation}
By Lemma \ref{tail} there exists $n_1$ so that for all $k$
\begin{equation}\label{kbound}
|Z_k(M)-\sum_{F \in \mathcal{S}(C)_k^{n_1+k}}
E(F)|<\frac{\epsilon}{4}.\end{equation}
Combining these, we get that for all $k>K$, all $n_0\geq n_1$

\begin{equation}|Z_{\infty}(M)-\sum_{F\in\mathcal{S}(C)_k^{k+n_0}}E(F)|<\frac{\epsilon}{2}.\end{equation}
Therefore, each of the $N-K-n_1-1$ terms of the sum
\begin{equation}\sum_{k=K+1}^{N-n_1}\sum_{F \in \mathcal{S}(C)_k^N} E(F)\end{equation}
is at most $\frac{\epsilon}{2}$ away from $Z_{\infty}(M)$. Since
$ \lim_{N\rightarrow \infty}\frac{N-K-n_1-1}{N}=1$ to finish the proof
it suffices to show that the first $K+1$ terms and the last $n_1$ terms inside the
outer sum on the right hand side of (\ref{regroup}) are
bounded regardless of the value of $N$.
For the first $K+1$ terms notice that  by (\ref{kbound})
\begin{equation}|\sum_{k=0}^K\sum_{F \in \mathcal{S}(C)_k^N} E(F)|<K(\frac{\epsilon}{4}+B),\end{equation}
where $B=\text{max}(|Z_0(M)|,|Z_1(M)|,\dots |Z_K(M)|)$.
The fact that the last $n_1$ inner sums 
\begin{equation}|\sum_{k=N-n_1}^N\sum_{F \in \mathcal{S}(C)_k^N} E(F)|\end{equation}
are bounded
regardless of $N$ follows from the fact that for every $i$
the limit 
\begin{equation}\lim_{k\rightarrow\infty}\sum_{F \in \mathcal{S}(C)_k^{k+i}} E(F)\end{equation}
exists (see Lemma \ref{tail}).

\item[{\bf(iii)}]
Absolute convergence of the sum (\ref{spherelim}) follows from the existence of the
universal bound on $|E_{\infty}(F)|$ for $F \in S_0(C)$. Since
$\lim_{k \rightarrow \infty}t^{4k}E_k(F)=E_{\infty}(F)$, this in turn
follows from a universal bound on $|t^{4k}E_k(F)|$ 
for $F \in S_0(C)$. 
By letting
$\epsilon=\frac{1}{2}$ in Remark
\ref{restatedlemma} we see that except for finitely many surfaces $F \in S_0(C)$,
$t^{4k}E_k(F)<\frac{1}{2}$. Since each of the sequences $t^{4k}E_k(F)$
is convergent for the remaining surfaces, the quantities $|t^{4k}E_k(F)|$ 
 are universally bounded for all surfaces $F \in S_0(C)$.

Our goal is to show that the sequence
\begin{equation}\label{finsum} 
 t^{8N}\sum_{F \in
    \mathcal{S}(C)^{2N}} E(F)
\end{equation} 
converges to  the sum (\ref{spherelim}). The first step is
to rewrite the 
finite sum in (\ref{finsum}) so that it is a sum over $0$-peelable
surfaces. We 
get,
\begin{equation} 
\sum_{F \in \cS_0(C)}\sum_{k=0}^{2N-m(F)} t^{8N}E_k(F).
\end{equation}
The largest part of this sum is at the end, so we change variables to put
the largest part at the beginning. Let $i=2N-m(F)-k$. Substitution,
along with splitting off an appropriate power of $t$, yields:
\begin{equation}\label{refinsum} 
\sum_{F \in \cS_0(C)}\sum_{i=0}^{2N-m(F)}t^{4m(F)+4i}
 t^{8N-4m(F)-4i}E_{2N-m(F)-i}(F).\end{equation}
From Remark \ref{restatedlemma} there exists $K_0$ so that, for all
$i\geq K_0$,
\begin{equation}  
\sum_{F \in \cS_0(C)^{T(K_0)}}t^{4i\chi(M)}E_i(F)<\frac{\epsilon (1-q)}{4}, 
\end{equation} thus
\begin{equation}  \sum_{F \in \cS_0(C)^{T(K_0)}}\frac{t^{4m(F)}}{1-q}
E_i(F)<\frac{\epsilon
  }{4}.\end{equation}
Estimating based on summing the geometric series $\sum_i t^{4i}=\frac{1}{1-q}$
we can truncate the sum (\ref{refinsum})  using any $K\geq K_0$  as follows and remain within 
$\epsilon/4$ of the
original sum.
\begin{equation}\label{firstrunc} 
\sum_{F \in \cS_0(C)^{K}}\sum_{i=0}^{2N-m(F)}t^{4m(F)+4i}
 t^{8N-4m(F)-4i}E_{2N-m(F)-i}(F).\end{equation}
Since by Proposition \ref{aug} the function $Q(F)$ is bounded below on
$\cS_0(C)$ there exists $B$
so that for all $F$, $N$ and $i$
\begin{equation} t^{8N-4m(F)-4i}E_{2N-m(F)-i}(F)<B.\end{equation}
From the elementary theory of the geometric series there exists
$I$ so that for all $F\in \cS_0(C)^{K}$,
\begin{equation} 
\sum_{i\geq I}^{2N-m(f)} t^{4m(F)+4i} B <\epsilon/4.\end{equation}
This means we can truncate the  sum (\ref{firstrunc}) again as follows and
remain within $\epsilon/4$ of 
the original sum:
\begin{equation}\label{secontruc}
 \sum_{F \in \cS_0(C)^{K}}\sum_{i=0}^{I}t^{4m(F)+4i}
 t^{8N-4m(F)-4i}E_{2N-m(F)-i}(F).\end{equation}
Using the fact that for any $F$, and for any fixed $i$,
\begin{equation} 
\lim_{N\rightarrow \infty}
t^{8N-4m(F)-4i}E_{2N-m(F)-i}(F)=E_{\infty}(F),
\end{equation}
together with the fact that the number of terms of the sum
(\ref{secontruc})  is bounded
independent of 
$N$, we can choose $N$ so large that the sum (\ref{secontruc}) is
within $\epsilon/4$ of
\begin{equation}\label{penultimate} 
\sum_{F \in \cS_0(C)^{K}}\sum_{i=0}^{I}t^{4m(F)+4i}
 E_{\infty}(F),\end{equation}
leaving us within $\frac{3\epsilon}{4}$ of the original sum
(\ref{finsum}). Using the 
absolute convergence of $\sum_{F \in \cS_0(C)} t^{4m(F)}E_{\infty}(F)$,
and the fact that the bound $B$ is still valid for $E_{\infty}(F)$,
we can choose $I$ large enough to make this last sum (\ref{penultimate}) 
within $\epsilon/4$ of 
\begin{equation} 
\sum_{F \in \cS_0(C)}\sum_{i=0}^{\infty}t^{4m(F)+4i}
 E_{\infty}(F).\end{equation}
Summing the geometric series yields the final result. \qed

\end{itemize}

\end{document}